\RequirePackage{fix-cm}
\documentclass[smallextended]{svjour3}       
\smartqed  

\usepackage[T1]{fontenc}
\usepackage[numbers,sort&compress]{natbib}
\usepackage{graphicx}
\usepackage{amsfonts,anysize,amsmath,indentfirst,geometry}
\usepackage{bm}
\usepackage{lineno}
\usepackage[colorlinks,citecolor=blue,linkcolor=blue]{hyperref}
\usepackage{mathrsfs}
\usepackage{amssymb}
\usepackage{array}
\usepackage{booktabs}
\usepackage{multirow}
\usepackage{algorithm}
\usepackage{algorithmic}
\usepackage{hyperref}
\usepackage{float}
\usepackage{subfigure}
\usepackage[justification=centering]{caption}
\usepackage{marvosym}
\usepackage{verbatim}
\usepackage{setspace}
\usepackage{lineno}

\numberwithin{equation}{section}
\newtheorem{thm}{Theorem}

\makeatletter
\newenvironment{breakablealgorithm}
  {
   \begin{center}
     \refstepcounter{algorithm}
     \hrule height.8pt depth0pt \kern2pt
     \renewcommand{\caption}[2][\relax]{
       {\raggedright\textbf{\ALG@name~\thealgorithm} ##2\par}%
       \ifx\relax##1\relax 
         \addcontentsline{loa}{algorithm}{\protect\numberline{\thealgorithm}##2}%
       \else 
         \addcontentsline{loa}{algorithm}{\protect\numberline{\thealgorithm}##1}%
       \fi
       \kern2pt\hrule\kern2pt
     }
  }{
     \kern2pt\hrule\relax
   \end{center}
  }
\makeatother

\DeclareGraphicsExtensions{.pdf,.jpeg,.jpg,.png,.eps}
\graphicspath{{./figures/},{./myfigures/}}
\begin{document}

\title{A rank-adaptive higher-order orthogonal iteration algorithm for truncated Tucker decomposition
}


\author{Chuanfu Xiao         \and
        Chao Yang 
}


\institute{Chuanfu Xiao \at
              School of Mathematical Sciences, Peking University, Beijing 100871, China \\
              \email{chuanfuxiao@pku.edu.cn}           
           \and
           Chao Yang (\Letter)\at
              School of Mathematical Sciences, Peking University, Beijing 100871, China\\
              National Engineering Laboratory for Big
Data Analysis \& Applications, Peking University, Beijing 100871, China\\
              \email{chao\_yang@pku.edu.cn}
}

\date{Oct 24, 2021}


\maketitle

\begin{abstract}
We propose a novel rank-adaptive higher-order orthogonal iteration (HOOI) algorithm to compute the truncated Tucker decomposition of higher-order tensors with a given error tolerance, and prove that the method is locally optimal and monotonically convergent.
A series of numerical experiments related to both synthetic and real-world tensors are carried out to show that the proposed rank-adaptive HOOI algorithm is advantageous in terms of both accuracy and efficiency.
Some further analysis on the HOOI algorithm and the classical alternating least squares method
are presented to further understand why rank adaptivity can be introduce into the HOOI algorithm and how it works.

\end{abstract}

\section{Introduction}
Tucker decomposition \cite{Tucker1966}, also known as the higher-order singular value decomposition (HOSVD) \cite{Lathauwer2000-1}, is regarded as a generalization of the classical singular value decomposition (SVD) for higher-order tensors. 
In practice, it usually suffices to consider the truncated Tucker decomposition, which can be seen as the low multilinear-rank approximation of higher-order tensors, as shown below:
\begin{equation}\label{prob:lowrankapp1}
\begin{gathered}
    \min\limits_{\bm{\mathcal{B}}}\|\bm{\mathcal{A}} - \bm{\mathcal{B}}\|_{F} \\
\mathrm{s.t.}\ \ \mu\mathrm{rank}(\bm{\mathcal{B}})\leq(R_{1},R_{2},\cdots,R_{N}),
\end{gathered}
\end{equation}
where $\bm{\mathcal{A}},\bm{\mathcal{B}}\in\mathbb{R}^{I_{1}\times I_{2}\cdots\times I_{N}}$ are $N$th-order tensors, $\|\cdot\|_{F}$ and $\mu\mathrm{rank}(\cdot)$ represent the Frobenius norm and the multilinear-rank of a higher-order tensor, respectively. Here in problem \eqref{prob:lowrankapp1},  the truncation $(R_{1},R_{2},\cdots,R_{N})$ is a suitably predetermined paramter. However, in many situations, the truncation is often very hard to obtain. Therefore, instead of problem \eqref{prob:lowrankapp1}, another form of the low multilinear-rank approximation problem is usually considered, which reads
\begin{equation}\label{prob:lowrankapp2}
    \begin{gathered}
    \min\limits_{\bm{\mathcal{B}}}\mu\mathrm{rank}(\bm{\mathcal{B}}) \\
            \mathrm{s.t.}\ \ \|\bm{\mathcal{A}} - \bm{\mathcal{B}}\|_{F}\leq\varepsilon\|\bm{\mathcal{A}}\|_{F},
    \end{gathered}
\end{equation}
where $\varepsilon$ is a given error tolerance. 

There are a plethora of algorithms to solve problem \eqref{prob:lowrankapp1}, such as the truncated HOSVD ($t$-HOSVD) algorithm \cite{Tucker1966,Lathauwer2000-1,Kolda2009}, the sequentially truncated HOSVD ($st$-HOSVD) algorithm \cite{Vannieuwenhoven2011,Vannieuwenhoven2012,Austin2016}, the higher-order orthogonal iteration (HOOI) algorithm \cite{Lathauwer2000-2,Kolda2009}, and iterative algorithms based on Riemannian manifold \cite{Elden2009,Ishteva2009,Savas2010,Ishteva2011}. But algorithms for solving problem \eqref{prob:lowrankapp2} are less to be seen; some representative 
work can be found in, e.g., Refs. \cite{Vannieuwenhoven2011,Vannieuwenhoven2012,Austin2016,ehrlacher2021adaptive}. 
In particular, Vannieuwenhove \emph{et al.} \cite{Vannieuwenhoven2011,Vannieuwenhoven2012} proposed to use uniform distribution strategies to determine the truncation based on the $t$-HOSVD and $st$-HOSVD algorithms. On top of the latter, Austin \emph{et al.} \cite{Austin2016} and Ballard \emph{et al.} \cite{ballard2020} studied high-performance implementation techniques for large-scale tensors on distributed-memory parallel computers. 
Recently, in order to further improve the accuracy of the truncation selected with the uniform distributed strategies, Ehrlacher \emph{et al.} \cite{ehrlacher2021adaptive} presented a greedy strategy for $t$-HOSVD, called Greedy-HOSVD, to search the truncation from a small initial truncation such as $(1,1,\cdots,1)$. 
Overall, the $t$-HOSVD and $st$-HOSVD algorithms have been playing a dominant role in the aforementioned references for solving problem \eqref{prob:lowrankapp2}, therefore these methods suffer from a same accuracy issue due to the fact that $t$-HOSVD and $st$-HOSVD are quasi-optimal algorithms for problem \eqref{prob:lowrankapp1} \cite{Vannieuwenhoven2012,Hackbusch2014,Minster2019}, often leading to inferior performance in applications. 


In this paper, we propose rank-adaptive HOOI, an accurate and efficient rank-adaptive algorithm for solving problem \eqref{prob:lowrankapp2}. Instead of relying on the $t$-HOSVD and $st$-HOSVD algorithms, we employ the HOOI method and propose a new strategy to adjust the truncation while updating the factor matrices for the truncated Tucker decomposition.
We prove that the proposed rank-adaptive HOOI algorithm satisfies local optimality, and is monotonically convergent, i.e., the truncation sequence is non-increasing during the iterations, and show by a series of experiments that it can outperform other algorithms in terms of both accuracy and efficiency.
We present further analysis on the HOOI algorithm to reveal that it is not an alternating least squares (ALS) method in the classical sense. 
Instead, the HOOI algorithm is essentially equivalent to the so called modified alternating least squares (MALS) method \cite{etter2016parallel,cichocki2016tensor} with orthogonal constraints, which is the main reason why 
rank adaptivity can be introduced.


The remainder of the paper is organized as follows. 
In Sec. \ref{sec:pre}, we introduce the HOOI algorithm, along with some basic notations for higher-order tensors. Then in Sec. \ref{sec:rank-ada hooi}, we present the rank-adaptive HOOI algorithm in detail, and prove that it is locally optimal and monotonically convergent. After that, a discussion on the HOOI algorithm is provided in Sec. \ref{sec:hooi}, revealing that HOOI is essentially equivalent to MALS with orthogonal constraints. Numerical experiments are reported in Sec. \ref{sec:experimets}, and the paper is concluded in Sec. \ref{sec:conslusions}.

\section{Overview of tensors and HOOI algorithm}\label{sec:pre}

In this paper, we use boldface capital calligraphic letters to represent higher-order tensors. 
Given an $N$th-order tensor $\bm{\mathcal{A}}\in\mathbb{R}^{I_{1}\times I_{2}\cdots\times I_{N}}$, the Frobenius norm is defined as
\[
\|\bm{\mathcal{A}}\|_{F} = \sqrt{\sum\limits_{i_{1},i_{2},\cdots,i_{N}}\bm{\mathcal{A}}_{i_{1},i_{2},\cdots,i_{N}}^{2}},
\]
where $\bm{\mathcal{A}}_{i_{1},i_{2},\cdots,i_{N}}$ denotes the $(i_{1},i_{2},\cdots,i_{N})$-th entry of the tensor. 
The mode-$n$ matricization of tensor $\bm{\mathcal{A}}$ is to reshape it to a matrix $\bm{A}_{(n)}\in\mathbb{R}^{I_{n}\times J_{n}}$, where $J_{n} = I_{1:N}/I_{n}$ and $I_{1:N} = \prod_{n=1}^{N}I_{n}$. Specifically, entry $\bm{\mathcal{A}}_{i_{1},i_{2},\cdots,i_{N}}$ of the tensor is mapped to the $(i_{n},j_{n})$-th entry of matrix $\bm{A}_{(n)}$, where
\[
j_{n} = 1+\sum\limits_{k=1,k\neq n}^{N}(i_{k}-1)J_{k}\ \ \mbox{with}\ \ J_{k} = \prod\limits_{m=1,m\neq n}^{k-1}I_{m}.
\]
The multiplication of a tensor $\bm{\mathcal{A}}$ and a matrix $\bm{U}\in\mathbb{R}^{J\times I_{n}}$ is denoted as $\bm{\mathcal{B}} = \bm{\mathcal{A}}\times_{n}\bm{U}$, which is also an $N$th-order tensor, i.e., $\bm{\mathcal{B}}\in\mathbb{R}^{I_{1}\cdots\times J\cdots\times I_{N}}$, and elementwisely we have
\[
\bm{\mathcal{B}}_{i_{1},\cdots,j,\cdots,i_{N}} = \sum\limits_{i_{n}=1}^{I_{n}}\bm{\mathcal{A}}_{i_{1},\cdots,i_{n},\cdots,i_{N}}\bm{U}_{j,i_{n}}.
\]
We denote the multilinear-rank of $\bm{\mathcal{A}}$ as $\mu\mathrm{rank}(\bm{\mathcal{A}})$, which is a positive integer tuple $(R_{1},R_{2},\cdots,R_{N})$, where $R_{n} = \mathrm{rank}(\bm{A}_{(n)})$ for all $n = 1,2,\cdots, N$. 

If the tensor $\bm{\mathcal{A}}$ has multilinear-rank $(R_{1},R_{2},\cdots,R_{N})$, then it can be expressed as
\begin{equation}\label{eq:tucker}
    \bm{\mathcal{A}} = \bm{\mathcal{G}}\times_{1}\bm{U}^{(1)}\times_{2}\bm{U}^{(2)}\cdots\times_{N}\bm{U}^{(N)},
\end{equation}
which is known as the Tucker decomposition of $\bm{\mathcal{A}}$, where $\bm{\mathcal{G}}\in\mathbb{R}^{R_{1}\times R_{2}\cdots\times R_{N}}$ is the core tensor, and $\bm{U}^{(n)}\in\mathbb{R}^{I_{n}\times R_{n}}$ ($n=1,2,\cdots,N$) are the column orthogonal factor matrices.
Based on the definition of Tucker decomposition \eqref{eq:tucker}, the low multilinear-rank approximation problem \eqref{prob:lowrankapp1} can be rewritten as
\begin{equation}\label{eq:low-rank tucker}
    \begin{gathered}
    \min\limits_{\bm{\mathcal{G}},\bm{U}^{(1)},\bm{U}^{(2)},\cdots,\bm{U}^{(N)}}\|\bm{\mathcal{A}} - \bm{\mathcal{G}}\times_{1}\bm{U}^{(1)}\times_{2}\bm{U}^{(2)}\cdots\times_{N}\bm{U}^{(N)}\|_{F}\\
            \mathrm{s.t.}\ \ (\bm{U}^{(n)})^{T}\bm{U}^{(n)}=\bm{I}\ \mathrm{for\ all}\ n=1,2,\cdots,N,
    \end{gathered}
\end{equation}
where $\bm{I}$ is a $R_{n}\times R_{n}$ identity matrix. 

As mentioned earlier, there are several algorithms to solve problem \eqref{eq:low-rank tucker}. 
Among them, the HOOI algorithm is the most popular iterative approach.
The computational precedure of HOOI is shown in Algorithm \ref{algo:hooi}.

\begin{breakablealgorithm}
	\algsetup{linenosize=\normalsize} 
	\setstretch{1.2}
	\normalsize
	\caption{The higher-order orthogonal iteration (HOOI) algorithm \cite{Lathauwer2000-2,Kolda2009}. }
	\label{algo:hooi}
	\begin{algorithmic}[1]
		\REQUIRE ~~\\
		Tensor $\bm{\mathcal{A}}\in\mathbb{R}^{I_{1}\times I_{2}\cdots\times I_{N}}$\\
		Truncation $(R_{1},R_{2},\cdots, R_{N})$ \\
		Initial guess $\{\bm{U}_{0}^{(n)}:n=1,2,\cdots,N\}$
		\ENSURE ~~\\
		Core tensor $\mathcal{G}$ \\
		Factor matrices $\{\bm{U}_{k}^{(n)}:n=1,2,\cdots,N\}$
		\STATE $k\ \leftarrow\ 0$
		\WHILE{not convergent} 
		\FORALL{$n\in\{1,2,\cdots,N\}$}
		\STATE $\bm{\mathcal{B}}\ \leftarrow\ \bm{\mathcal{A}}\times_{1}(\bm{U}_{k+1}^{(1)})^{T}\cdots\times_{n-1}(\bm{U}_{k+1}^{(n-1)})^{T}\times_{n+1}(\bm{U}_{k}^{(n+1)})^{T}\cdots\times_{N}(\bm{U}_{k}^{(N)})^{T}$
		\STATE $\bm{B}_{(n)}\ \leftarrow\ \bm{\mathcal{B}}$ in matrix format 
		\STATE $\bm{U},\bm{\Sigma},\bm{V}^{T}\ \leftarrow$ truncated rank-$R_{n}$ SVD of $\bm{B}_{(n)}$
		\STATE $\bm{U}_{k+1}^{(n)}\ \leftarrow \bm{U}$
		\STATE $k\ \leftarrow\ k+1$
		\ENDFOR
		\ENDWHILE
		\STATE $\bm{\mathcal{G}}\ \leftarrow\ \bm{\Sigma}\bm{V}^{T}$ in tensor format
	\end{algorithmic}
\end{breakablealgorithm}

It is easy to see that HOOI is an alternating iterative algorithm, in which the low multilinear-rank approximation of $\bm{\mathcal{A}}$ is calculated by alternately updating the factor matrices $\{\bm{U}^{(n)}:n=1,2,\cdots,N\}$ through rank-$R_{n}$ approximations of the matricized tensors $\{\bm{B}_{(n)}:n=1,2,\cdots,N\}$. 

\section{Rank-adaptive HOOI algorithm}\label{sec:rank-ada hooi}

Analogous to \eqref{eq:low-rank tucker}, the low multilinear-rank approximation problem \eqref{prob:lowrankapp2} can be reformulated as:
\begin{equation}\label{eq:low-rank approx-2}
    \begin{array}{r}
    \min\limits_{\scriptsize
    \begin{array}{c}
    \bm{\mathcal{G}},\bm{U}^{(1)},\bm{U}^{(2)},\cdots,\bm{U}^{(N)}\\
    \|\bm{\mathcal{A}} - \bm{\mathcal{G}}\times_{1}\bm{U}^{(1)}\times_{2}\bm{U}^{(2)}\cdots\times_{N}\bm{U}^{(N)}\|_{F}\leq\varepsilon\|\bm{\mathcal{A}}\|_{F}
    \end{array}}(R_{1},R_{2},\cdots,R_{N})\\
            \mathrm{s.t.}\ \  (\bm{U}^{(n)})^{T}\bm{U}^{(n)} = \bm{I}\ \mathrm{for\ all}\ n=1,2,\cdots,N.
    \end{array}
\end{equation}

Existing approaches for solving problem \eqref{eq:low-rank approx-2} are all based on the $st$-HOSVD algorithms, which are quasi-optimal \cite{Vannieuwenhoven2012,Hackbusch2014,Minster2019}
in the sense that 
\[
    \begin{gathered}
        \|\bm{\mathcal{A}} - \bm{\mathcal{A}}_{\mathrm{opt}}\|_{F}\leq\|\bm{\mathcal{A}} - \bm{\mathcal{A}}_{(s)t-\mathrm{HOSVD}}\|_{F}\leq\sqrt{N}\|\bm{\mathcal{A}} - \bm{\mathcal{A}}_{\mathrm{opt}}\|_{F},
    \end{gathered}
\]
where $\bm{\mathcal{A}_{\mathrm{opt}}}$ represents the best low multilinear-rank approximation of $\bm{\mathcal{A}}$, and $\bm{\mathcal{A}}_{(s)t-\mathrm{HOSVD}}$ is the low multilinear-rank approximation obtained by the $t$- or $st$-HOSVD algorithm. Therefore the truncation obtained with existing methods is usually much larger than the exact one, often resulting in inferior performance in applications. 

To address this issue, we tackle problem \eqref{eq:low-rank approx-2} from a different angle based on the HOOI algorithm. 
On top of Algorithm \ref{algo:hooi}, we propose a new rank adaptive strategy to adjust the truncation.
Specifically, during the HOOI iterations, $R_{n}$ for dimension $n$ is updated by minimizing $R$ that satisfies
\begin{equation}\label{eq:select truncation}
\|\bm{B}_{(n)} - (\bm{B}_{(n)})_{[R]}\|_{F}^{2}\leq\|\bm{\mathcal{B}}\|_{F}^{2} - (1-\varepsilon^{2})\|\bm{\mathcal{A}}\|_{F}^{2},
\end{equation} 
where $(\bm{B}_{(n)})_{[R]}$ is the best rank-$R$ approximation of $\bm{B}_{(n)}$. 
This can be done by calculating the full matrix SVD of $\bm{B}_{(n)}$.
The update of factor matrix $\bm{U}^{(n)}$ inherits the spirits of the original HOOI algorithm, which is based on the leading $R_{n}$ left singular vectors of $\bm{B}_{(n)}$. The detailed procedure of the rank-adaptive HOOI algorithm is presented in Algorithm \ref{algo:rank-ada hooi}.


\begin{breakablealgorithm}
	\algsetup{linenosize=\normalsize} 
	\setstretch{1.2}
	\normalsize
	\caption{Rank-adaptive HOOI algorithm}
	\label{algo:rank-ada hooi}
	\begin{algorithmic}[1]
		\REQUIRE ~~\\
		Tensor $\bm{\mathcal{A}}\in\mathbb{R}^{I_{1}\times I_{2}\cdots\times I_{N}}$\\
		Error tolerance $\varepsilon$ \\
		Initial truncation $(R_{1}^{0},R_{2}^{0},\cdots,R_{N}^{0})$ \\
		Initial guess $\{\bm{U}_{0}^{(n)}:n=1,2,\cdots,N\}$ 
		\ENSURE ~~\\
		Truncation $(R_{1}^{k},R_{2}^{k},\cdots,R_{N}^{k})$ \\
		Core tensor $\mathcal{G}_k$ \\
		Factor matrices $\{\bm{U}_{k}^{(n)}:n=1,2,\cdots,N\}$
   	    \STATE $\bm{\mathcal{G}}_{0}\ \leftarrow\  \bm{\mathcal{A}}\times_{1}(\bm{U}_{0}^{(1)})^{T}\times_{2}(\bm{U}_{0}^{(2)})^{T}\cdots\times_{N}(\bm{U}_{0}^{(N)})^{T}$
		\STATE $k\ \leftarrow\ 0$
		\WHILE{$\|\bm{\mathcal{G}}_{k}\|_{F} > \sqrt{1-\varepsilon}\|\bm{\mathcal{A}}\|_{F}$}
		\FORALL{$n\in\{1,2,\cdots,N\}$}
		\STATE $\bm{\mathcal{B}}\ \leftarrow\ \bm{\mathcal{A}}\times_{1}(\bm{U}_{k+1}^{(1)})^{T}\cdots\times_{n-1}(\bm{U}_{k+1}^{(n-1)})^{T}\times_{n+1}(\bm{U}_{k}^{(n+1)})^{T}\cdots\times_{N}(\bm{U}_{k}^{(N)})^{T}$
		\STATE $\bm{B}_{(n)}\ \leftarrow\ \bm{\mathcal{B}}$ in matrix format
		\STATE $\bm{U},\bm{\Sigma},\bm{V}\ \leftarrow$ full-SVD of $\bm{B}_{(n)}$
		\STATE $R_{n}^{k}\ \leftarrow$ {minimize} $R$ such that $\sum\limits_{r>R}\bm{\Sigma}_{r,r}^{2}\leq\|\bm{\mathcal{B}}\|_{F}^{2} - (1-\varepsilon^{2})\|\bm{\mathcal{A}}\|_{F}^{2}$
		\STATE $\bm{U}^{(n)}_{k+1}\ \leftarrow\ \bm{U}_{:,1:R^{k}_{n}}$
		\ENDFOR
		\STATE $\bm{\mathcal{G}}_{k+1}\ \leftarrow\ \bm{\Sigma}_{1:R^{k}_{n},1:R^{k}_{n}}\bm{V}_{:,1:R_{n}^{k}}^{T}$ in tensor format
		\STATE $k\ \leftarrow\ k+1$
		\ENDWHILE
	\end{algorithmic}
\end{breakablealgorithm}

We remark that for Algorithm \ref{algo:rank-ada hooi}, the initial guess $\{\bm{U}_{0}^{(1)},\bm{U}_{0}^{(2)},\cdots,\bm{U}_{0}^{(N)}\}$ needs to be a feasible solution to problem \eqref{prob:lowrankapp2}.
This can be achieved with low cost by using the $st$-HOSVD algorithm or randomized methods such as the ones suggested in Refs. \cite{che2019randomized,Minster2019,ahmadi2021randomized}.



The reason for using Eq. \eqref{eq:select truncation} to update truncation $R_{n}$ is twofold. Firstly, Eq. \eqref{eq:select truncation} ensures the low multilinear-rank approximation $\bm{\mathcal{G}}_{k}\times_{1}\bm{U}_{k}^{(1)}\times_{2}\bm{U}_{k}^{(2)}\cdots\times_{N}\bm{U}_{k}^{(N)}$ satisfies 
\[
    \|\bm{\mathcal{A}} - \bm{\mathcal{G}}_{k}\times_{1}\bm{U}_{k}^{(1)}\times_{2}\bm{U}_{k}^{(2)}\cdots\times_{N}\bm{U}_{k}^{(N)}\|_{F}\leq\varepsilon\|\bm{\mathcal{A}}\|_{F}
\]
for all $k$, which means $\bm{\mathcal{G}}_{k}\times_{1}\bm{U}_{k}^{(1)}\times_{2}\bm{U}_{k}^{(2)}\cdots\times_{N}\bm{U}_{k}^{(N)}$ is a feasible solution of problem \eqref{prob:lowrankapp2}. Secondly and more importantly, Eq. \eqref{eq:select truncation} is a local optimal strategy for updating $R_{n}$, 
which is illustrated by Theorem \ref{thm:local-optimal} as follows.

\begin{thm}\label{thm:local-optimal}
	Let $(R_{1}^{k+1},\cdots,R^{k+1}_{n-1},R_{n}^{k},\cdots,R_{N}^{k})$ be the truncation before updating $R_{n}$, and the corresponding factor matrices are $\{\bm{U}^{(1)}_{k+1},\cdots,\bm{U}_{k+1}^{(n-1)},\bm{U}_{k}^{(n)},\cdots,\bm{U}_{k}^{(N)}\}$, then 
	\begin{equation}\label{eq:sub-problem}
    \begin{split}
    &R_{n}^{k+1}=\mathop{\arg\min}\limits_{R}\|\bm{B}_{(n)} - (\bm{B}_{(n)})_{[R]}\|_{F}^{2}\leq\|\bm{\mathcal{B}}\|_{F}^{2} - (1-\varepsilon^{2})\|\bm{\mathcal{A}}\|_{F}^{2}\\
\end{split}
\end{equation}
is the optimal strategy for updating $R_{n}$, where $\bm{B}_{(n)}$ is the mode-$n$ matricization of 
\[
\bm{\mathcal{B}} = \bm{\mathcal{A}}\times_{1}(\bm{U}_{k+1}^{(1)})^{T}\cdots\times_{n-1}(\bm{U}_{k+1}^{(n-1)})^{T}\times_{n+1}(\bm{U}_{k}^{(n+1)})^{T}\cdots\times_{N}(\bm{U}_{k}^{(N)})^{T}.
\]
\end{thm}

\noindent\textbf{Proof.} Given factor matrices $\{\bm{U}^{(1)}_{k+1},\cdots,\bm{U}_{k+1}^{(n-1)},\bm{U}_{k}^{(n+1)},\cdots,\bm{U}_{k}^{(N)}\}$, using the expression of Tucker decomposition \eqref{eq:tucker}, problem \eqref{prob:lowrankapp2} can be rewritten as
\begin{equation}\label{eq:local-optimal-1}
    \begin{gathered}
    \min\mathrm{rank}(\bm{G}_{(n)})\\
            \mathrm{s.t.}\ \ \|\bm{\mathcal{A} - \bm{\mathcal{G}}}\times_{1}\bm{U}_{k+1}^{(1)}\cdots\times_{n-1}\bm{U}_{k+1}^{(n-1)}\times_{n+1}\bm{U}_{k}^{(n+1)}\cdots\times_{N}\bm{U}_{k}^{(N)}\|_{F}\leq\varepsilon\|\bm{\mathcal{A}}\|_{F}.
    \end{gathered}
\end{equation}
Because the factor matrices are column orthogonal and the Frobenius norm satisfies the orthogonal invariance property, the constraint condition of \eqref{eq:local-optimal-1} is equivalent to
\[
\begin{split}
    \|\bm{\mathcal{B}} - \bm{\mathcal{G}}\|_{F}^{2}+\|\bm{\mathcal{A}} - \bm{\mathcal{B}}\times_{1}\bm{U}_{k+1}^{(1)}\cdots\times_{n-1}\bm{U}_{k+1}^{(n-1)}\times_{n+1}\bm{U}_{k}^{(n+1)}\cdots\times_{N}\bm{U}_{k}^{(N)}\|_{F}^{2}
    \leq\varepsilon^{2}\|\bm{\mathcal{A}}\|_{F}^{2},
    \end{split}
\]
which further leads to
\begin{equation}\label{eq:local-optimal-2}
\begin{split}
    \|\bm{\mathcal{A}}\|_{F}^{2}-\|\bm{\mathcal{B}}\|_{F}^{2}+\|\bm{\mathcal{B}} - \bm{\mathcal{G}}\|_{F}^{2}
    \leq\varepsilon^{2}\|\bm{\mathcal{A}}\|_{F}^{2}.
    \end{split}
\end{equation}
Therefore, Eq. \eqref{eq:local-optimal-1} can be reformulated as 
\[\label{eq:local-optimal-3}
    \begin{gathered}
    \min\mathrm{rank}(\bm{G}_{(n)})\\
            \mathrm{s.t.}\ \ \|\bm{B}_{(n)} - \bm{G}_{(n)}\|_{F}^{2}\leq\|\bm{\mathcal{B}}\|_{F}^{2} - (1-\varepsilon^{2})\|\bm{\mathcal{A}}\|_{F}^{2}.
    \end{gathered}
\]
Since $R_{n}^{k+1}=\mathop{\arg\min}\limits_{R}\|\bm{B}_{(n)} - (\bm{B}_{(n)})_{[R]}\|_{F}^{2}\leq\|\bm{\mathcal{B}}\|_{F}^{2} - (1-\varepsilon^{2})\|\bm{\mathcal{A}}\|_{F}^{2}$, for any positive integer $R<R_{n}^{k+1}$, we have
\[
    \|\bm{B}_{(n)} - (\bm{B}_{(n)})_{[R]}\|_{F}^{2}>\|\bm{\mathcal{B}}\|_{F}^{2} - (1-\varepsilon^{2})\|\bm{\mathcal{A}}\|_{F}^{2},
\]
where $(\bm{B}_{(n)})_{[R]}$ is the best rank-$R$ approximation of $\bm{B}_{(n)}$. This means that for any core tensor $\bm{G}_{(n)}$ with rank less than $R_{n}^{k+1}$, we have 
\[
\|\bm{B}_{(n)} - \bm{G}_{(n)}\|_{F}^{2}>\|\bm{\mathcal{B}}\|_{F}^{2} - (1-\varepsilon^{2})\|\bm{\mathcal{A}}\|_{F}^{2},
\]
which implies that $R_{n}^{k+1}$ is the optimal solution of sub-problem \eqref{eq:local-optimal-1}. $\hfill\square$


From Theorem \ref{thm:local-optimal}, it is known that \eqref{eq:sub-problem} is the optimal strategy for updating $R_{n}$. However, the update of the factor matrix $\bm{U}^{(n)}$ is not unique. Different strategies to update $\bm{U}^{(n)}$ could lead to different algorithm behaviors.  Theorem \ref{thm:non-increasing} illustrates that the truncation sequence $\{\bm{R}^{k}\}$ is non-increasing as iteration proceeds, if the factor matrix $\bm{U}^{(n)}$ is obtained by the leading $R_{n}^{k+1}$ left singular vectors of $\bm{B}_{(n)}$, where $\bm{R}^{k}=(R_{1}^{k},R_{2}^{k},\cdots,R_{N}^{k})$ is the truncation in the $k$-th iteration.

\begin{thm}\label{thm:non-increasing}
	Let $\{\bm{R}^{k}\}$ be the truncation sequence during the iterations of Algorithm \ref{algo:rank-ada hooi}, then 
	\begin{equation}\label{eq:rank non-increasing}
	\bm{R}^{k+1}\leq\bm{R}^{k},
	\end{equation}
	i.e., $R_{n}^{k+1}\leq R_{n}^{k}$ for all $n=1,2,\cdots,N$, and the infimum of $\{\bm{R}^{k}\}$ can be reached. 
\end{thm}


\noindent\textbf{Proof.} Let $\{\bm{U}^{(1)}_{k+1},\cdots,\bm{U}_{k+1}^{(n-1)},\bm{U}_{k}^{(n)},\cdots,\bm{U}_{k}^{(N)}\}$ and $\bm{\mathcal{B}}$ be defined as in Theorem \ref{thm:local-optimal}, and $\bm{U}\in\mathbb{R}^{I_{n}\times R_{n}^{k}}$ be composed of the leading $R_{n}^{k}$ left singular vectors of $\bm{B}_{(n)}$, we have 
\begin{equation}\label{eq: thm non-increasing 1}
\begin{gathered}
		\|\bm{\mathcal{A}} - \bm{\mathcal{G}}\times_{1}\bm{U}^{(1)}_{k+1}\cdots\times_{n-1}\bm{U}^{(n-1)}_{k+1}\times_{n}\bm{U}\cdots\times_{N}\bm{U}^{(N)}_{k}\|_{F}\leq\\
		 \|\bm{\mathcal{A}} - \bm{\mathcal{G}}_{n-1,k}\times_{1}\bm{U}^{(1)}_{k+1}\cdots\times_{n-1}\bm{U}^{(n-1)}_{k+1}\times_{n}\bm{U}^{(n)}_{k}\cdots\times_{N}\bm{U}^{(N)}_{k}\|_{F}\leq\varepsilon\|\bm{\mathcal{A}}\|_{F},
	\end{gathered}
\end{equation}
where $\bm{\mathcal{G}}_{n-1,k} = \bm{\mathcal{B}}\times_{n}(\bm{U}^{(n)}_{k})^{T}$ and $\bm{\mathcal{G}} = \bm{\mathcal{B}}\times_{n}\bm{U}^{T}$. 
Reformulating Eq. \eqref{eq: thm non-increasing 1}, we have
\[
    \|\bm{B}_{(n)} - \bm{U}\bm{U}^{T}\bm{B}_{(n)}\|_{F}^{2}\leq\|\bm{\mathcal{B}}\|_{F}^{2}-(1-\varepsilon^{2})\|\bm{\mathcal{A}}\|_{F}^{2}.
\]
Since $\bm{U}\in\mathbb{R}^{I_{n}\times R_{n}^{k}}$ is composed of the leading $R_{n}^{k}$ left singular vectors of $\bm{B}_{(n)}$, the rank of  $\bm{U}\bm{U}^{T}\bm{B}_{(n)}$ must be equal to $R_{n}^{k}$. According to \eqref{eq:sub-problem}, it is easy to know that $R_{n}^{k+1}\leq R_{n}^{k}$. Further, because $R_{n}^{k}$ is a positive integer for all $n$ and $k$, the infimum of $\{\bm{R}^{k}\}$ can be reached. $\hfill\square$



\section{Discussion of HOOI algorithm}\label{sec:hooi}


Besides HOOI, the classical ALS method can also be used to solve problem \eqref{prob:lowrankapp1}, in which the core tensor and the factor matrices $\{\bm{\mathcal{G}},\bm{U}^{(1)},\bm{U}^{(2)},\cdots,\bm{U}^{(N)}\}$ are updated sequentially in a certain order, such as 
\[
    \bm{\mathcal{G}}\rightarrow\bm{U}^{(1)}\rightarrow\bm{U}^{(2)}\cdots\rightarrow\bm{U}^{(N)}.
\]
When updating a factor matrix $\bm{U}^{(n)}$, the core tensor and the other factor matrices are fixed, and the update of the core tensor $\bm{\mathcal{G}}$ is done by fixing all factor matrices.
Overall, the core tensor and the factor matrices are updated one after another along the ALS iteration.

However, HOOI and the above classical ALS method are often referred to with different names and even occasionally mistaken with each other.
Table \ref{table:name} shows the different names of them in some related references.
For example, 
the TUCKER-TS proposed in Ref. \cite{malik2018low} is in fact the randomization of the classical ALS method, instead of HOOI as claimed in the paper, via a tensor sketch technique.
And the P-Tucker in Ref. \cite{oh2018scalable} and the general Tucker factorization algorithm (GTA) in Ref. \cite{Oh2019} are the essentially parallel implementations of the classical ALS method, rather than HOOI as shown in the papers, on CPU and GPU platforms, respectively. 
On the other hand, the TUCKALS3 algorithm presented in Ref. \cite{Kroonenberg1980} is equivalent to the HOOI algorithm for third-order tensors, 
which was further extended to $N$th-order tensors and called TuckALS$n$ in Ref. \cite{kapteyn1986approach}.
Later, the HOOI algorithm was officially proposed in Ref. \cite{Lathauwer2000-2} and was recognized as an efficient iterative approach to update the factor matrices \cite{Kolda2009}. 
Nowadays, the HOOI algorithm is sometimes referred to as Tucker-ALS; examples can be found in Refs. \cite{Oh2017,ma2018accelerating,ma2021fast}.

\begin{table}[H]
	\normalsize
	\begin{center}
		\caption{\normalsize Different names of classical ALS and HOOI in some references.}\label{table:name}
		\begin{tabular}{c|c|c|c|c|c|c}
			\toprule
			\multicolumn{1}{c}{} &\multicolumn{3}{|c}{Classical ALS} & \multicolumn{3}{|c}{HOOI} \\
			\midrule
			Name & TUCKER-TS & P-Tucker & GTA & TuckALS$n$ & HOOI & Tucker-ALS  \\
			\midrule
			Ref.  & \cite{malik2018low}  & \cite{oh2018scalable} & \cite{Oh2019} & \cite{Kroonenberg1980,kapteyn1986approach} & \cite{Lathauwer2000-2,Kolda2009}  & \cite{Oh2017,ma2018accelerating,ma2021fast}   \\
			\bottomrule
		\end{tabular}
	\end{center}
\end{table}

Here we would like to point out that, though quite alike, the HOOI algorithm is essentially different from the classical ALS method in the way that the sub-problem is defined and solved. In particular, Theorem \ref{thm:subproblem} illustrates that the sub-problem of the classical ALS method is an ordinary least squares problem with multiple right-hand sides, while the sub-problem of HOOI is a constrained least squares problem.

\begin{thm}\label{thm:subproblem}
	In the classical ALS algorithm, the sub-problem for updating the factor matrix $\bm{U}^{(n)}$ is an ordinary least squares problem with multiple right-hand sides: 
	\begin{equation}\label{eq:als-subpro}
	    \begin{gathered}
	        \min\limits_{\bm{X}\in\mathbb{R}^{I_{n}\times R_{n}}}\|\bm{A}_{(n)}-\bm{X}\bm{M}_{1}^{T}\|_{F}^{2},
	    \end{gathered}
	\end{equation}
	where 
	$\bm{M}_{1} = \bm{U}^{(N)}\cdots\otimes\bm{U}^{(n+1)}\otimes\bm{U}^{(n-1)}\cdots\otimes\bm{U}^{(1)}\bm{G}_{(n)}^{T}$. 
	And in the HOOI algorithm, the sub-problem is a constrained least squares problem:
	\begin{equation}\label{eq:subproblem-hooi}
    \begin{gathered}
    \min\limits_{\bm{X}}\|\bm{A}_{(n)} - \bm{X}\bm{M}_{2}^{T}\|_{F}^{2},\quad
            \mathrm{s.t.}\ \ \mathrm{rank}(\bm{X})\leq R_{n},
    \end{gathered}
\end{equation}
where 
$\bm{M}_{2} = \bm{U}^{(N)}\cdots\otimes\bm{U}^{(n+1)}\otimes\bm{U}^{(n-1)}\cdots\otimes\bm{U}^{(1)}$. 
\end{thm}

\noindent\textbf{Proof.} In the classical ALS method, when the core tensor and the factor matrices other than $\bm{U}^{(n)}$ are fixed, it is easy to know, from the multilinear property of the Tucker decomposition, that the low multilinear rank approximation problem degenerates to the ordinary least squares problem \eqref{eq:als-subpro}.
In the HOOI method, the factor matrices are column orthogonal during the iteration, then the Kronecker product of them is also column orthogonal, which means that the coefficient matrix of \eqref{eq:subproblem-hooi} is column orthogonal. Suppose that $\bm{Q} = \left[\bm{M}_{2},\bm{M}_{2}^{\bot}\right]$ is an orthogonal matrix, due to the orthogonal invariance of the Frobenius norm, we have 
\[
    \begin{gathered}
        \|\bm{A}_{(n)} - \bm{X}\bm{M}_{2}^{T}\|_{F}^{2} = \|(\bm{A}_{(n)} - \bm{X}\bm{M}_{2}^{T})\bm{Q}\|_{F}^{2}
        = \|\bm{A}_{(n)}\bm{M}_{2} - \bm{X}\|_{F}^{2} + \|\bm{A}_{(n)}\bm{M}_{2}^{\bot}\|_{F}^{2},
    \end{gathered}
\]
where the $\|\bm{A}_{(n)}\bm{M}_{2}^{\bot}\|_{F}^{2}$ is not dependent on $\bm{A}_{(n)}$. Thus the problem \eqref{eq:subproblem-hooi} is equivalent to the following rank-$R_{n}$ approximation problem,
\begin{equation}\label{eq:low-rank}
    \begin{aligned}
    \min\limits_{\bm{X}}\|\bm{X}-\bm{A}_{(n)}\bm{M}_{2}\|_{F}^{2} \quad
        \begin{array}{r@{\quad}l@{}l@{\quad}l}
            \mathrm{s.t.} & \mathrm{rank}(\bm{X})\leq R_{n}.
        \end{array}
    \end{aligned}
\end{equation}
And since $\bm{A}_{(n)}\bm{M}_{2}$ is the mode-$n$ matricization of $\bm{\mathcal{B}}$, where
\[
\bm{\mathcal{B}} = \bm{\mathcal{A}}\times_{1}(\bm{U}^{(1)})^{T}\cdots\times_{n-1}(\bm{U}^{(n-1)})^{T}\times_{n+1}(\bm{U}^{(n+1)})^{T}\cdots\times_{N}(\bm{U}^{(N)})^{T},
\]
then solving problem \eqref{eq:low-rank} is equivalent to solving the following problem:
\[
    \begin{gathered}
        \max\limits_{\bm{U}^{(n)}\in\mathbb{R}^{I_{n}\times R_{n}}}\|\bm{\mathcal{A}}\times_{1}(\bm{U}^{(1)})^{T}\times_{2}(\bm{U}^{(2)})^{T}\cdots\times_{N}(\bm{U}^{(N)})^{T}\|_{F}^{2}\\
        \mathrm{s.t.}\ \ (\bm{U}^{(n)})^{T}\bm{U}^{(n)} = \bm{I},
    \end{gathered}
\]
which is the sub-problem of HOOI. Therefore, the constrained least squares problem \eqref{eq:subproblem-hooi} is the sub-problem of HOOI. $\hfill\square$

It is clear that the sub-problem \eqref{eq:als-subpro} is equivalent to a series of ordinary least squares problems, which can be solved in parallel with a row-wise rule \cite{oh2018scalable,Oh2019}. For the sub-problem \eqref{eq:subproblem-hooi}, Theorem \ref{thm:hooi-mals} illustrates that it can be solved by truncated matrix SVD.

\begin{thm}\label{thm:hooi-mals}
	In the HOOI algorithm, the constrained least squares problem \eqref{eq:subproblem-hooi} can be solved by calculating the truncated rank-$R_{n}$ SVD of $\bm{A}_{(n)}\bm{M}_{2}$, i.e., $\bm{B}_{(n)}$. In this way, the core tensor $\bm{\mathcal{G}}$ and the factor matrix $\bm{U}^{(n)}$ are updated simultaneously during the iterations. 
\end{thm}
\noindent\textbf{Proof.} Since the sub-problem of HOOI is the rank-$R_{n}$ approximation of $\bm{B}_{(n)}$, and by the Ecart-Young theorem \cite{eckart1936approximation}, we know that the rank-$R_{n}$ SVD of $\bm{B}_{(n)}$ is a solution of it. Suppose that $\bm{U},\bm{\Sigma},\bm{V}$ are the truncated rank-$R_{n}$ SVD factors, then $\bm{U}$ is the update value of $\bm{U}^{(n)}$, and
\begin{equation}\label{eq:thmhooi-mals1}
    \bm{\Sigma}\bm{V}^{T} = (\bm{U}^{(n)})^{T}\bm{B}_{(n)}.
\end{equation}
It is easy to know that the right-hand side of Eq. \eqref{eq:thmhooi-mals1} is the mode-$n$ matricization of $\bm{\mathcal{B}}\times_{n}(\bm{U}^{(n)})^{T}$, that is, 
\[
\bm{\mathcal{A}}\times_{1}(\bm{U}^{(1)})^{T}\times_{2}(\bm{U}^{(2)})^{T}\cdots\times_{N}(\bm{U}^{(N)})^{T},
\]
which is the updated core tensor $\bm{\mathcal{G}}$. Therefore, the core tensor $\bm{\mathcal{G}}$ and the factor matrix $\bm{U}^{(n)}$ are updated simultaneously by solving the sub-problem \eqref{eq:subproblem-hooi}. $\hfill\square$

Theorem \ref{thm:hooi-mals} shows that the HOOI algorithm updates the core tensor with a factor matrix simultaneously by solving a constrained least squares problem in each iteration, which can be calculated by truncated matrix SVD.
And this results indicate that it is possible for 
rank adaptation based on the intermediate tensor $\bm{\mathcal{B}}$. Specifically, one can adaptively adjust the truncation according to Eq. \eqref{eq:sub-problem} to solve the low multilinear-rank approximation problem with a given error tolerance, i.e., problem \eqref{prob:lowrankapp2}. On the other hand, it is difficult to introduce any certain of adaptivity into the classical ALS method. We believe that, due largely to the fact that the HOOI algorithm was often confused with classical ALS in practical applications, the rank-adaptive HOOI algorithm has not surfaced until today.

We further remark that, the update rule for updating the two factors simultaneously is sometimes called modified ALS (MALS) in tensor computation references \cite{etter2016parallel,cichocki2016tensor,cichocki2017tensor}, which is equivalent to the 
two-site density matrix renormalization group (DMRG) algorithm in the field of quantum physics \cite{schollwock2011density,legeza2014tensor,cichocki2016tensor}. Clearly, Theorem \ref{thm:hooi-mals} also illustrates that HOOI is essentially equivalent to an MALS method with orthogonal constraints, which require that the factor matrices are column orthogonal during the HOOI iteration.
This property of orthogonal constraints is indispensable for successful rank adaptation; 
otherwise, the coefficient matrix $\bm{M}_{2}$ in the constrained least squares problem \eqref{eq:subproblem-hooi} will not be column orthogonal, thus the solution is not the rank-$R_{n}$ approximation of $\bm{B}_{(n)}$. However, in the classical ALS method, the orthogonal constraints are not necessarily satisfied, which is the essential difference from HOOI.
In addition, the core tensor is updated more frequently in the HOOI method, so intuitively, HOOI converges faster than classical ALS.

\section{Numerical experiments}\label{sec:experimets}

In this section, we will examine the performance of the proposed rank-adaptive HOOI algorithm with several numerical experiments related to both synthetic and real-world tensors.
For the purpose of comparison, we choose the $t$-HOSVD and $st$-HOSVD algorithms provided in Tensorlab 3.0 \cite{Vervliet} with a uniform distribution strategy to determine the truncation, and implement the Greedy-HOSVD algorithm with a bottom-up greedy strategy to search the truncation from $(1,1,\cdots,1)$ \cite{ehrlacher2021adaptive}.
We also complement Greedy-HOSVD with a new truncation selection strategy to start the truncation searching from the multilinear-rank of original input tensor in a top-down manner.
The implementations of the Greedy-HOSVD algorithms and our proposed rank-adaptive HOOI are all based on Tensorlab 3.0 to make fair comparison.
When using the proposed rank-adaptive HOOI algorithm, we set the maximum number of iterations to $500$ and compute the initial guess using both $st$-HOSVD and a randomized method with randomly generated orthogonal matrices \cite{che2019randomized,Minster2019,ahmadi2021randomized}.
All the experiments are carried out on a computer equipped with an Intel Xeon Gold 6240 CPU of 2.60 GHz and MATLAB R2019b.


\subsection{Reconstruction of a low multilinear-rank tensor with Gaussian noise} 
The first test case is to reconstruct a low multilinear-rank tensor with Gaussian noise. We generate the input tensor via
\[
\bm{\hat{\mathcal{A}}} = \bm{\mathcal{A}}/\|\bm{\mathcal{A}}\|_{F}+\delta\bm{\mathcal{E}}/\|\bm{\mathcal{E}}\|_{F},
\]
where $\bm{\mathcal{A}}\in\mathbb{R}^{500\times500\times500}$ is a tensor with $\mu\mathrm{rank}(\bm{\mathcal{A}}) = (100,100,100)$, the elements of $\bm{\mathcal{E}}$ follow the standard Gaussian distribution, and the noise level is set to $\delta = 10^{-2}$. 
We run the test using algorithms including $t$-HOSVD, $st$-HOSVD, Greedy-HOSVD with both bottom-up and top-down adaptive strategies, and rank-adaptive HOOI initialized with both $st$-HOSVD and randomization.
The error tolerance for tensor reconstruction
is set to $\epsilon = 10^{-2}$,
where
$\epsilon = \|\bm{\mathcal{B}} - \bm{\mathcal{A}}/\|\bm{\mathcal{A}}\|_{F}\|_{F}$, 
and $\bm{\mathcal{B}}$ is the low multilinear-rank approximation of $\bm{\hat{\mathcal{A}}}$ obtained by the tested algorithms.

\begin{table}[h!]
	\begin{center}
		\caption{Reconstruction error, truncation and running time of various algorithms for reconstructing a Gaussian noisy low multilinear-rank tensor.}\label{table:random data}
		\begin{tabular}{c|c|c|c|c}
			\toprule
		\multicolumn{2}{c|}{Algorithms} &  Reconstruction error   & Truncation & Running time (s) \\
			\midrule
	\multicolumn{2}{c|}{$t$-HOSVD} & $8.0535\times10^{-3}$ & (441,485,369) & 8.45 \\
	\multicolumn{2}{c|}{$st$-HOSVD} & $7.8057\times10^{-3}$ & (441,484,346) & 7.81 \\
			\midrule
			 \multirow{2}{*}{Greedy-HOSVD} & bottom-up & $5.4748\times10^{-3}$ & (327,327,329) & 8.09 \\
			 & top-down & $5.4505\times10^{-3}$ & (326,326,328) & 8.06 \\
			\midrule
			\multirow{2}{*}{Rank-adaptive HOOI}  & $st$-HOSVD & $0.9466\times10^{-3}$ & (100,100,100) & 12.66 \\
			& random & $0.9466\times10^{-3}$ & (100,100,100) & 6.34 \\
			\bottomrule
		\end{tabular}
	\end{center}
\end{table}

Table \ref{table:random data} shows the reconstruction error, truncation and running time of the tested algorithms. 
From the table, it can be seen that 
the performances of the bottom-up and top-down strategies for the Greedy-HOSVD are very close with each other, both failing to reach the exact truncation, though better than plain $t$- and $st$-HOSVD.
And the rank-adaptive HOOI algorithm can successfully find the exact truncation no matter how the initial guess is calculated, while the truncations obtained by other algorithms are far from the exact one.
As a result, the reconstruction error obtained with the rank-adaptive HOOI method is much smaller than that of other tested algorithms. 
The time cost of rank-adaptive HOOI
is also competitive, especially when the initial guess is obtained randomly.


\subsection{Compression of a regularized Coulomb kernel} 

The second test case is to compress a fourth-order tensor, which is constructed from a function of four variables  $f:\mathbb{R}^{4}\rightarrow\mathbb{R}$ defined as \cite{ehrlacher2021adaptive}:
\[
    f(x_{1},x_{2},x_{3},x_{4}) = \ln{(0.1+\sqrt{|x_{1}-x_{2}|^{2}}+\sqrt{|x_{3}-x_{4}|^{2}})}.
\]
The function can be seen as a regularization of the 2D Coulomb kernel. In this example, the variable $x_{i}$ takes $I$ values in the interval $[-100,100]$ for all $i=1,2,3,4$, thus the input tensor is a fourth-order tensor that belongs to $\mathbb{R}^{I\times I\times I\times I}$. 
We test algorithms including $t$-HOSVD, $st$-HOSVD, greedy-HOSVD with both bottom-up and top-down adaptive strategies, and rank-adaptive HOOI initialized with both $st$-HOSVD and randomization.


\begin{table}[h!]
	\begin{center}
		\caption{Number of parameters required by various algorithms for compressing the regularized Coulomb kernel with different compression levels.}\label{table:tolerance}
		\scalebox{1.0}{\begin{tabular}{c|c|c|c|c|c|c}
			\toprule
		\multirow{2}{*}{$\varepsilon$} &  \multirow{2}{*}{$t$-HOSVD} &  \multirow{2}{*}{$st$-HOSVD} & 
		\multicolumn{2}{|c}{Greedy-HOSVD} &\multicolumn{2}{|c}{Rank-adaptation HOOI}\\
			\cline{4-7}
	& &  &	bottom-up & top-down & $st$-HOSVD &  random  \\
			\midrule
			1.0e-1 & 5128 & 2928 & 2708 & 1824 & 1616  & 1616 \\
			1.0e-2 & $2.05\times10^{7}$ & $2.61\times10^{6}$ & $2.21\times10^{5}$ & $1.84\times10^{5}$ &  $3.37\times10^{4}$ & $5.20\times10^{4}$\\
			1.0e-3 & $1.48\times10^{9}$ & $1.44\times10^{9}$ & $1.32\times10^{9}$ & $1.30\times10^{9}$  & $1.01\times10^{9}$ & $1.02\times10^{9}$ \\
			\bottomrule
		\end{tabular}}
	\end{center}
\end{table}

\begin{figure}[h!]
	\centering
	\includegraphics[width=1.0\hsize]{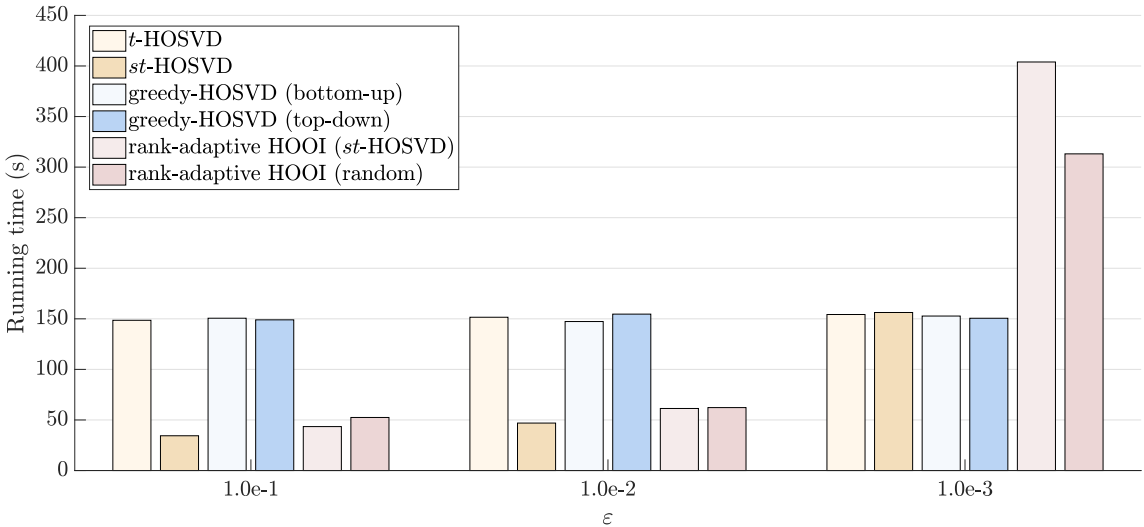}
	\caption{Running time of various algorithms for compressing the regularized Coulomb kernel with different compression levels.}\label{fig:coulomb1}
\end{figure}

First, we fix the dimension size $I=200$ and compress the tensor under three different error tolerances, namely $\varepsilon = 10^{-1},10^{-2},10^{-3}$, corresponding to high, medium, and low compression levels, respectively.
The required number of parameters as well as the running time of the tested algorithm are shown in Table \ref{table:tolerance} and Figure \ref{fig:coulomb1}, respectively. 
From Table \ref{table:tolerance}, we observe that Greedy-HOSVD leads to
compressing parameters much less than $t$- and $st$-HOSVD algorithms do, especially when the top-down strategy is applied, while the proposed rank-adaptive HOOI algorithm can further effectively reduce the number of parameters.
Figure \ref{fig:coulomb1} illustrates that the running time of proposed rank-adaptive HOOI algorithm is competitive as compared with other choices, except for the case of low compression level with $\varepsilon = 10^{-3}$.
We remark that the compression level with $\varepsilon = 10^{-3}$ is in fact too low to be of any practical value because in this case the numbers of parameters after compression are close to the size of full tensor, which is $1.6\times10^{9}$.


\begin{table}[h!]
	\begin{center}
		\caption{Number of parameters required by various algorithms for compressing the regularized Coulomb kernel with different tensor sizes.}\label{table:size}
		\scalebox{1.0}{\begin{tabular}{c|c|c|c|c|c|c}
			\toprule
		\multirow{2}{*}{$I$} &  \multirow{2}{*}{$t$-HOSVD} &  \multirow{2}{*}{$st$-HOSVD} & 
		\multicolumn{2}{|c}{Greedy-HOSVD} &\multicolumn{2}{|c}{Rank-adaptation HOOI}\\
			\cline{4-7}
	& &  &	bottom-up & top-down & $st$-HOSVD &  random  \\
			\midrule
			100 & $3.03\times10^{7}$ & $1.95\times10^{7}$ & $3.20\times10^{6}$ & $2.91\times10^{6}$ & $1.22\times10^{6}$  & $1.40\times10^{6}$ \\
			200 & $2.05\times10^{7}$ & $2.61\times10^{6}$ & $2.21\times10^{5}$ & $1.84\times10^{5}$ & $3.37\times10^{4}$  & $5.20\times10^{4}$\\
			300 & $8.83\times10^{6}$ & $9.87\times10^{5}$ & $1.60\times10^{5}$ & $1.32\times10^{5}$ & $4.12\times10^{4}$ & $4.12\times10^{4}$ \\
			\bottomrule
		\end{tabular}}
	\end{center}
\end{table}

\begin{figure}[h!]
	\centering
	\includegraphics[width=1.0\hsize]{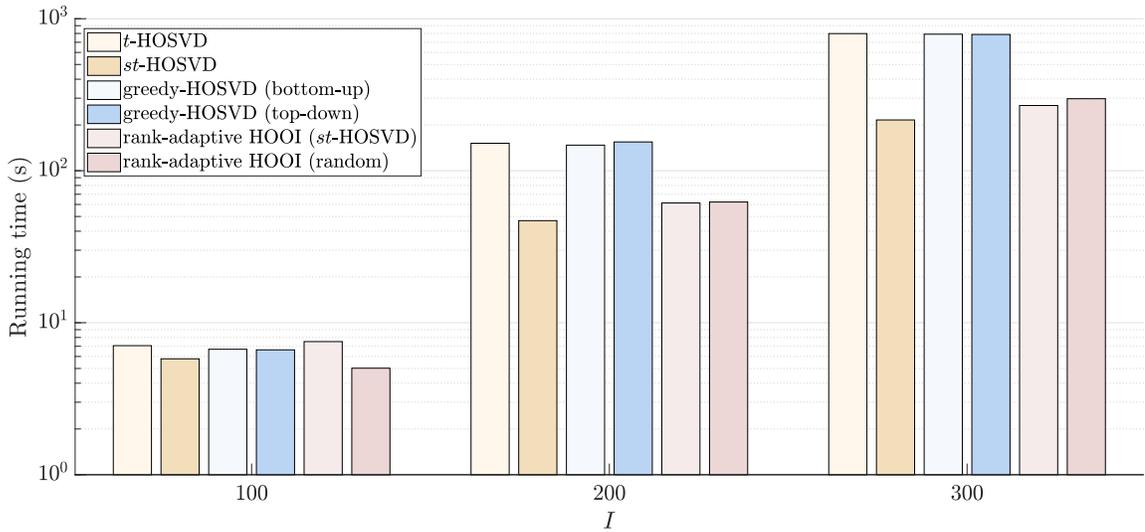}
	\caption{Running time of various algorithms for compressing the regularized Coulomb kernel with different tensor sizes.}\label{fig:coulomb2}
\end{figure}

Then, we fix the error tolerance $\varepsilon=10^{-2}$ and adjust the tensor size to be $I=100, 200, 300$ to run the test again. The number of parameters and running time of the tested algorithms are shown in Table \ref{table:size} and Figure \ref{fig:coulomb2}, respectively.
From Table \ref{table:size}, we find that as the tensor size $I$ grows, the rank-adaptive HOOI algorithm maintains a clear advantage in the number of parameters compared to other methods. For example, when rank-adaptive HOOI is initialized with $st$-HOSVD, the number of parameters can be been reduced by $2.4\times\sim24.8\times$, $5.5\times\sim60.8\times$ and $3.2\times\sim214.3\times$ for $I=100, 200, 300$, respectively. From Figure \ref{fig:coulomb2}, we observe again that the time cost of rank-adaptive HOOI is kept among low in all tested algorithms.

\subsection{Classification of handwritten digits data}

The third test case is the classification of handwritten digits data, and the input tensor is from the MNIST database \cite{MNIST}. We use a fourth-order tensor $\bm{\mathcal{A}}\in\mathbb{R}^{28\times28\times5000\times10}$ to represent the training dataset, where the first and second modes of $\bm{\mathcal{A}}$ are the texel modes, the third mode corresponds to training images, and the fourth mode represents image categories. We employ methods introduced in Ref. \cite{Savas2007} to compress the training data of images by Tucker decomposition and utilize the core tensor to classify on the test data. 
We set the error tolerance for the low multilinear-rank approximation to 
$\epsilon = 0.45$ and
run the test using algorithms
including $t$-HOSVD, $st$-HOSVD, Greedy-HOSVD with both bottom-up and top-down adaptive strategies, and rank-adaptive HOOI initialized with both $st$-HOSVD and randomization. 


\begin{figure}[h!]
	\centering
	\includegraphics[width=1.0\hsize]{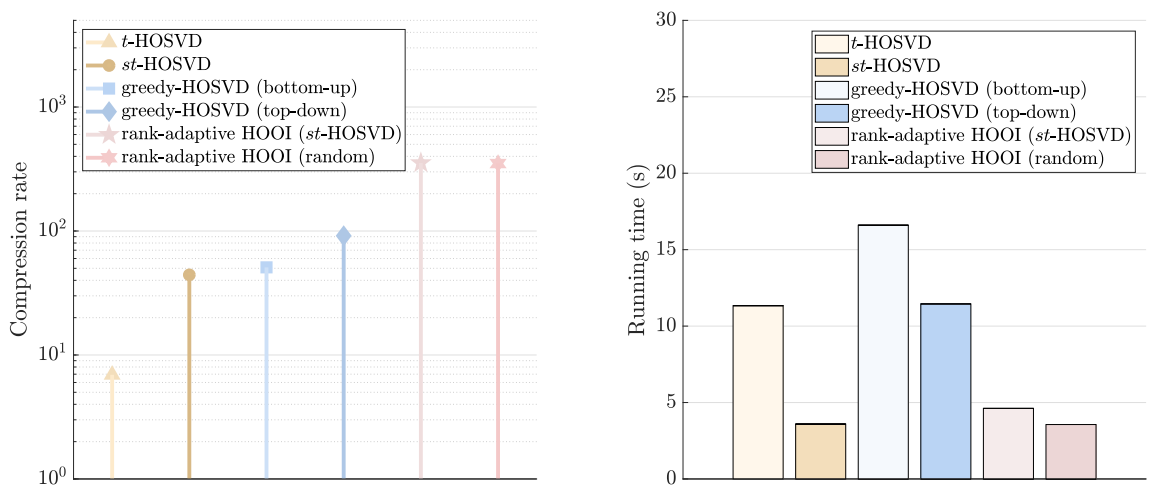}
	\caption{Compression rate and running time of various algorithms for compressing handwritten digits data.}\label{fig:mnist} 
\end{figure}

To examine the compression effect, we calculate the compression rate that is defined as $\mathrm{size}(\bm{\mathcal{A}})/\mathrm{size}(\bm{\mathcal{G}})=\prod_{n=1}^{N}I_{n}/R_{n}$.
Figure \ref{fig:mnist} shows the compression rate and running time of the tested algorithms.
It can be observed from the figure that Greedy-HOSVD has better compression rate than the classical $t$- and $st$-HOSVD algorithms, especially when applying the top-down adaptive strategy, and the proposed rank-adaptive HOOI algorithm is able to lead to a further significant improvement. Specifically, the compression rate obtained with rank-adaptive HOOI is $12.11\times\sim51.17\times$ higher than $t$-HOSVD, $1.89\times\sim8.00\times$ higher than $st$-HOSVD, and $1.65\times\sim6.96\times$ higher than Greedy-HOSVD. 

\begin{table}[h!]
	\begin{center}
		\caption{Accuracy and running time of various algorithms for classifying handwritten digits data.}\label{table:mnist}
		\scalebox{1.0}{
		\begin{tabular}{c|c|c|c}
			\toprule
		\multicolumn{2}{c|}{Algorithms} &  Classification accuracy ($\%$)  & Classification time (s)  \\
			\midrule
	\multicolumn{2}{c|}{$t$-HOSVD} & 93.40 & 44.85 \\
	\multicolumn{2}{c|}{$st$-HOSVD} & 94.99 &  20.73 \\
			\midrule
			 \multirow{2}{*}{Greedy-HOSVD} & bottom-up & 94.11 &  101.51 \\
			 & top-down & 94.89 &  71.77 \\
			\midrule
			\multirow{2}{*}{Rank-adaptive HOOI}  & $st$-HOSVD & 93.21 &  0.41 \\
			& random & 93.03 &  0.37 \\
			\bottomrule
		\end{tabular}}
	\end{center}
\end{table}

We then examine the performance of classification based on the compressed tensors obtained by the tested algorithms.
The results on the accuracy and running time for classification are shown in Table \ref{table:mnist}.
It is seen that the classification accuracy is close when using different algorithms, staying between $93\%$ to $95\%$. 
Because of the higher compression rate that leads to smaller core tensor for classification, the rank-adaptive HOOI algorithm has a much faster classification speed. This result clearly shows the advantage of the proposed rank-adaptive HOOI algorithm over other methods.

\section{Conclusions}\label{sec:conslusions}

In this paper, a new rank-adaptive HOOI algorithm is presented to efficiently compute the truncated Tucker decomposition of higher-order tensors with a given error tolerance. 
We prove the local optimality and monotonic convergence of the proposed rank-adaptive HOOI method and show by a series of numerical experiments the advantages of it.
Further analysis on the HOOI algorithm is also provided to reveal the difference between HOOI and clasical ALS and understand why rank adaptivity can be introduced.
As future work, we plan to conduct in-depth convergence analysis on the rank-adaptive HOOI algorithm and study the corresponding parallel computing techniques on modern high-performance computers.


\section*{Declarations}
\noindent\textbf{Funding} This study was funded in part by 
 the National Key Research and Development Program of China (\#2018AAA0103304).
\smallskip

\noindent\textbf{Conflicts of Interest} The authors declare that they have no conflict of interest.
\smallskip

\noindent\textbf{Availability of Data and Material} The datasets generated and analyzed during the current study are available from 
the corresponding author on reasonable request.
\smallskip

\noindent\textbf{Code Availability} The code used in the current study is available from the corresponding author on reasonable request.
\smallskip


%
%

\bibliographystyle{spmpsci}      
\bibliography{article.bib}   

%
%

\end{document}